\documentclass[12pt]{article} 
\usepackage{amssymb} 
\usepackage{amsthm} 
\usepackage{amsgen} 
\usepackage{amsfonts} 
\usepackage{amscd}
 
\newtheorem{theorem}{Theorem}[section]  
\newtheorem{corollary}[theorem]{Corollary}  
  
\newtheorem{question}[theorem]{Question}  
\newtheorem{lemma}[theorem]{Lemma}  
\newtheorem{proposition}[theorem]{Proposition}

\newtheorem{remark}[theorem]{Remark}  
\def\fff{\mathbb{F}}

\def\ppp{\mathbb{P}}

\def\pf{{\bf proof}:\ }
\def\qed{$\Box$}

\begin{document} 

\author{David Joyner and Pablo Lejarraga} 

\title{On the variety of Borels in relative position $\vec{w}$} 
\date{1-31-2003} 
\maketitle 

\begin{abstract} 
Let $G$ be a connected semi-simple group defined over 
and algebraically closed field, $T$ a fixed Cartan,
$B$ a fixed Borel containing $T$,
$S$ a set of simple reflections associated to the
simple positive roots corresponding to $(T,B)$,
and let ${\cal B}\cong G/B$ denote the Borel variety.
For any $s_i\in S$, $1\leq i\leq n$, let
$\overline{O}(s_1,\dots ,s_n)= 
\{(B_0,\dots ,B_{n})\in {\cal B}^{n+1}\ |\ 
(B_{i-1},B_{i})\in \overline{O(s_i)},\ 1\leq i\leq n\}$,
where $O(s)$ denotes the subvariety of pairs of Borels in
${\cal B}^2$ in relative position $s$. 
We show that such varieties are smooth and indicate why this result is,
in one sense, best possible.
Our main results assume that $k$ has characteristic $0$.
\end{abstract} 

\tableofcontents 

\newpage

Let $G$ be a connected semi-simple group defined over 
and algebraically closed field, $T$ a fixed Cartan,
$B$ a fixed Borel containing $T$,
$S$ a set of simple reflections associated to the
simple positive roots corresponding to $(T,B)$,
$W$ the Weyl group of $T$,
and let ${\cal B}\cong G/B$ denote the Borel variety.
For $w_i\in W$, let $\vec{w}=(w_1,\dots ,w_n)$
and $\overline{O}(\vec{w})= 
\{(B_0,\dots ,B_{n})\in {\cal B}^{n+1}\ |\ 
(B_{i-1},B_{i})\in \overline{O(w_i)},\ 1\leq i\leq n\}$,
where $O(w)$ denotes the subvariety of pairs of Borels in
${\cal B}^2$ in relative position $w$ (defined below)
and $\overline{O(w)}$
denotes Zariski closure of ${O(w)}$
in ${\cal B}\times {\cal B}$. 
We call such a variety the 
{\bf closed variety of Borels of relative position $\vec{w}$}.

Such varieties have appeared in many
works. We mention, for example, \cite{Han},
\cite{DL}, and \cite{Spr}. In special cases, they are
known to be non-singular. 
Let $\vec{s}=(s_1,\dots ,s_n)\in S^n$.
We prove that these varieties $\overline{O}(\vec{s})$ are
non-singular.
Although this seems to be a simple result, we did not succeed
in finding a simple proof. 
The problem lies in the distinction between 
reduced fibers and schematic fibers in the statement of
a theorem of Mumford (Theorem 3' in \S III.10 of \cite{Mum})
which we use at one point in the proof.
Most of the technicalities, which ultimately rely on results of
Grothendieck \cite{EGA}, have been pushed into
\S \ref{sec:appendix}.

We also prove an analogous smoothness result for the corresponding 
Deligne-Lusztig varieties, under a technical assumption. 
In some cases, such varieties are known to have a large number of rational
points, so there is the hope that such results will be useful for constructing
error-correcting codes having ``good'' parameters \cite{Han2}.

Finally, we indicate why our smoothness result is best possible in 
some sense. Consider the projection 
$\pi_w: \overline{O}(w)\rightarrow \overline{O}(\emptyset)={\cal B}$,
where $w\in W$ is {\it not} a simple reflection.
In general, it is known that $\overline{O}(w)$ is singular 
(this is the raison d'etre for the Demazure-Hansen 
resolution \cite{DL}, \S 9.1)
and of course $\overline{O}(\emptyset)$ (the Borel variety) is smooth.
Therefore, by Lemma \ref{lemma:1.1} below, $\pi_w$ cannot
be smooth. As a consequence of this type of reasoning,
our smoothness result cannot be generalized to
the analogous projection
$\overline{O}(w_1,\dots ,w_n)
\rightarrow \overline{O}(w_1,\dots ,w_{n-1})$,
$w_i\in W$.

{\it Ackowledgement}: We thank S. Hansen for useful 
correspondence on \cite{Han}, \cite{Han2}.

\section{Main results} 
\label{sec:1}

Let $G$ be a connected semisimple 
group defined over an algebraically closed field 
$k$ 
with maximal torus $T$ contained in a Borel $B$. 
This is given the structure of a Coxeter group
generated by the reflections $s_\alpha$ associated to
the simple positive roots $\alpha\in \Delta$ of $(T,B)$. 
We identify $G$ with the group of $k$-valued points of $G$. 

We shall denote conjugation of group elements by $\ ^g x=gxg^{-1}$. 
Representing $w\in W$ by
$\dot{w} \in N_G(T)$, define 
\[
^w B' =\dot{w} B' {\dot{w}}^{-1},
\]
for any Borel $B'$ of $G$.
The Borel subgroups belonging to the set
$\{\ ^w B\ |\ w\in W\}$ are precisely the Borel
subgroups of $G$ which contain $T$. 

Let ${\cal B}\cong G/B $ denote the Borel variety. 

We say that
$B',B''\in {\cal B}$ are in {\bf position $w$} if
$(B',B'')$ belongs to the $G$-orbit of $(B , \, ^w B )$, for $w\in W$. 
We denote this orbit by 
\[
O(w)= 
\{ (B',B'')\in {\cal B}^2\ |\ B',B''\ {\rm are\ in \ position}\ w\}. 
\] 
In particular, $O(1)\cong {\cal B}$ is the diagonal in 
${\cal B}^2$. This orbit $O(w)$ is open in its Zariski closure 
$\overline{O(w)}\subset {\cal B}\times {\cal B}$ (this follows
from a much more general result in \cite{St}, using the
fact that ${\cal B}$ has a transitive $G$-action). 
It is known that 
\[ 
\overline{O(w)}=\cup_{w'\leq w}O(w'), 
\] 
where $\leq$ denotes the Bruhat ordering 
associated to $(T ,B)$ (\cite{Spr}, \S 10.2.13). 

Let $\vec{w}=(w_1,\dots ,w_n)\in W^n$ 
and let 
\[ 
O(\vec{w})
= \{(B_0,\dots ,B_{n})\in {\cal B}^{n+1}\ |\ 
(B_{i-1},B_{i})\in {O(w_i)},\ 1\leq i\leq n\}. 
\] 
We call this the
{\bf variety of Borels of relative position $\vec{w}$}.

We say an algebraic group $G$ {\bf operates} on a variety $X$
if there is a morphism $\phi:G\times X\rightarrow X$
satisfying $\phi(1,x)=x$ for all $x\in X$
(which may or may be compatible with the group 
operation). We say $G$ {\bf acts} on $X$ if it operates
on $X$ in a way compatible with the group operation.
For notational simplicity, we usually suppress the $\phi$ from the
notation and use juxtaposition to denote an operation.

A basic fact about this variety is the following.

\begin{lemma}
\label{lemma:action}
Let $\vec{s}=(s_1,\dots ,s_n)\in S^n$.
The variety $O(\vec{s})$ is the orbit of a connected 
algebraic group operation. 
\end{lemma}

This is proven in \S \ref{sec:2}.

Let 
\[ 
\overline{O}(w_1,\dots ,w_n)= 
\{(B_0,\dots ,B_{n})\in {\cal B}^{n+1}\ |\ 
(B_{i-1},B_{i})\in \overline{O(w_i)},\ 1\leq i\leq n\}. 
\] 
This variety is closed in ${\cal B}^{n+1}$. 
In Lemma \ref{lemma:proj}, we will see that this may be written 
as an iterated fiber product of the $\overline{O}(w_i)$'s..
The projection morphism 
\[ 
\pi_n:\overline{O}(w_1,\dots ,w_n)\rightarrow \overline{O}(w_1,\dots ,w_{n-1}), 
\] 
is defined\footnote{
If $n=1$, we interpret 
$\overline{O}(w_1,\dots ,w_{0})=\overline{O}(\emptyset)$ 
as ${\cal B}$.} 
by ignoring the last ``coordinate''.
We show in \S \ref{sec:3} that $\pi_n$ is flat.
However, in the case when the $w_i$ are simple reflections, much more is
true.

\begin{proposition}
\label{prop:smooth} 
Assume that $k$ has characteristic $0$.
Let $\vec{s}=(s_1,\dots ,s_n)\in S^n$.
For $n\geq 1$, $\overline{O}(s_1,\dots ,s_n)$ is 
non-singular, and $\pi_n$ is smooth. 
\end{proposition}

This is proven in \S \ref{sec:3}.
We show in \S \ref{sec:2} and \S \ref{sec:appendix} that each 
schematic fiber of $\pi_n$ is $\cong \ppp^1$
(for ${\rm char}(k)=0$ - in fact, this result is needed to prove the
above proposition and is the reason for the characteristic $0$ hypothesis
there). We believe that these results hold even in 
characteritic $p>0$.

\begin{corollary}
If $w = s_1\dots s_n$ is a minimal expression then
the projection map 
$(B_0,\dots ,B_{n})\longmapsto (B_0,B_{n})$ induces a 
surjection 
the morphism $\pi : \overline{O}(s_1,\dots ,s_n)\rightarrow \overline{O(w)}$ 
such that 
\begin{equation} 
\pi:\pi^{-1}(O(w))\rightarrow O(w) 
\label{eq:1} 
\end{equation} 
is an isomorphism.
In other words, $\pi$ is 
a resolution of singularities of $\overline{O(w)}$.

\end{corollary}

This is well-known (see, for example, \S 9 in Deligne-Lusztig \cite{DL}
or \S 10.2 in Springer \cite{Spr}). 
If, in addition, the $s_1$, \dots  ,$s_n$ are distinct then
this resolution is an isomorphism \cite{Han}.

\begin{question} Is
\label{thrm:main} 
\[
\overline{O}(s_1,\dots ,s_n)\stackrel{\pi_n}{\rightarrow}
\overline{O}(s_1,\dots ,s_{n-1})
\]
a $\ppp^1$-bundle?
\end{question} 

Some remarks on this are given in \S \ref{sec:4}.

Let $\vec{s}=(s_1,\dots ,s_n)\in S^n$ be a list of
simple reflections, with respect to
a fixed Borel $B$ in a connected 
semisimple group (functor) $G$ defined over
a finite field $\fff_q$. Let 
$F:G\rightarrow G$ denote the Frobenius morphism.
Following Hansen \cite{Han}, let 
\[
\begin{array}{c}
\overline{X}(\vec{s})
=\\
\{ (\, ^{g_0}B,\dots ,\, ^{g_n}B)\in {\cal B}^{n+1}\ |\ 
g_k^{-1}g_{k+1}\in B\cup Bs_{k+1}B,\ {\rm for\ }0\leq k<n,
g_n\in F(g_0)\}.
\end{array}\] 
Here we shall call 
$\overline{X}(\vec{s})$ a {\bf Deligne-Lusztig variety}.

\begin{corollary}
\label{cor:DL}
Assume that Proposition \ref{prop:smooth} 
holds for $k=\overline{\fff_q}$.
$\overline{X}(\vec{s})$ is non-singular.
\end{corollary}

This is proven in \S \ref{sec:5}.

\section{Definitions and lemmas}
\label{sec:2}

We begin with some simple observations on terminology.
We call $\overline{O}(s_1,\dots ,s_n)$
the closed variety of Borels of relative position $\vec{w}$.
On the other hand, the
closure of the variety of Borels of relative position $\vec{w}$
is $\overline{O(s_1,\dots ,s_n)}$. How are these related?

\begin{lemma}
$\overline{O}(s_1,\dots ,s_n)=\overline{O(s_1,\dots ,s_n)}$.
\end{lemma}

\begin{remark}
It is clear that 
$\overline{O}(s_1,\dots ,s_n)$ is closed in ${\cal B}^{n+1}$
since it is defined by closed conditions.
This implies 
$\overline{O(s_1,\dots ,s_n)}
\subset \overline{O}(s_1,\dots ,s_n)$.
By Lemma \ref{lemma:action} above and 
the proof of Lemma 2.3.3(i) in \cite{Spr2},
$O(s_1,\dots ,s_n)$ contains a set which is open in 
$\overline{O(s_1,\dots ,s_n)}$. 
Observe that $O(s_1,\dots ,s_n)$ is 
open in $\overline{O}(s_1,\dots ,s_n)$, since its
complement is closed.
Thus, both $\overline{O}(s_1,\dots ,s_n)$ and
$\overline{O(s_1,\dots ,s_n)}$
contain $O(s_1,\dots ,s_n)$ as a dense subset. 

\end{remark}

\pf
For notational simplicity, let $X={\cal B}$, let
$Y_i=O(s_i)\subset X\times X$, and let
$f_i:X^{n+1}\rightarrow Y_i$ be the 
morphism $f_i(x_0,\dots ,x_n)=(x_{i-1},x_i)$,
for $1\leq i\leq n$. Let
\[
Z=\{x\in X^{n+1}\ |\ f_i(x)\in Y_i,\forall 1\leq i\leq n\}.
\]
To prove the lemma, it suffices to prove
\begin{equation}
\overline{Z}=\{x\in X^{n+1}\ |\ f_i(x)\in \overline{Y_i},
\forall 1\leq i\leq n\}.
\label{eqn:Zbar}
\end{equation}
Let $g:X^{n+1}\rightarrow Y_1\times \dots \times Y_n$,
so $g=(f_1,\dots ,f_n)$. Since each $f_i$ is proper
(and hence closed), so is $g$. Then
$Z=g^{-1}(Y_1\times \dots \times Y_n)$.
We have
\[
g^{-1}(Y_1\times \dots \times Y_n)
\subset \overline{g^{-1}(Y_1\times \dots \times Y_n)}
\subset
g^{-1}(\overline{Y_1\times \dots \times Y_n}),
\]
since $g$ is continuous in the Zariski topology.
Recall $G$ is connected, so each $Y_i=O(s_i)$
is irreducible. Since $g$ is also closed, 
it follows that $g^{-1}(Y_1\times \dots \times Y_n)$
is open in its closure. 
This implies (\ref{eqn:Zbar}).
\qed

Next we determine the reduced fibers of $\pi_n$. 

\begin{proposition}
\label{lemma:fiber}
For $n\geq 1$, the projection morphism 
induces a surjective morphism 
\[ 
\pi_n:\overline{O}(s_1,\dots ,s_n)\rightarrow 
\overline{O}(s_1,\dots ,s_{n-1}), 
\] 
whose reduced fibers are isomorphic to $\ppp ^1$. 
\end{proposition}

\begin{remark}
\label{remark:irred}
(1) The morphism $\pi_n$ is proper, since it can be expressed as the
composition of a closed immersion (inclusion) with a base-change morphism
(projection ${\cal B}^{n+1}\rightarrow {\cal B}^n$).

(2) In the case when $w=s_1\dots s_n$ is a minimal expression, see
for example, \cite{Spr}, \S 10.2.13.
\end{remark}

\pf
We want to show the projection map
\[
\overline{O}(s_1,\dots ,s_n) \rightarrow \overline{O}(s_1,\dots ,s_{n-1})
\]
is surjective with fibers $\ppp^1$. 
Observe that there is an isomorphism of reduced fibers
\begin{equation}
\label{eqn:fibers}
\pi_n^{-1}(B_0,\dots ,B_{n-1})
\rightarrow 
\pi_1^{-1}(B_{n-1}),
\end{equation}
given by 
$(B'_0,...,B'_{n})\longmapsto (B'_{n-1},B'_n)$.
This boils the proof down to investigating
\[
\pi_1^{-1}(B') = (B'\times {\cal B})\cap \overline{O}(s)=
\{(B' ,B'')\ |\ (B' ,B'') \in \overline{O(s)}\},
\]
where $B' $ is fixed and $s$ is simple.

For $\alpha$ a root of $B'$, let $U_\alpha$ denote the
unipotent subgroup of $B' $ associated to $\alpha$
and $T$ the maximal torus of $B'$.
Note 
\[
B' = T  \times \prod_{\alpha >0} U_{\alpha},
\]
and 
\[
^sB' = T \times \prod_{s(\alpha) >0} U_{\alpha},
\]
where $>$ denotes the Bruhat ordering with respect to $(B',T)$.
Let $P_s=\langle B' ,\, ^sB' \rangle $ be the parabolic
generated by $B'$ and $^sB'$. 
For $(B' ,B'')\in \pi_1^{-1}(B')$, we have $B''=pB' p^{-1}$, for some
$p\in P_s$ (the
conjugating element must belong to the
parabolic generated by the two Borels by Corollary 11.17(iii)
in \cite{Bor}). Let $Y=P_s/B' $.
By \cite{St}, \S 3.9, we have $Y\cong \ppp^1$.
Since a Borel is its own normalizer, we have
an isomorphism $Y\rightarrow \pi_1^{-1}(B')$.
The result follows.
\qed

\begin{corollary} 
\label{cor:irred}
$\overline{O}(s_1,\dots ,s_n)$ is irreducible.
\end{corollary}

\pf
We assume for the moment that $\pi_n$ is flat
(this will be proven in Lemma \ref{lemma:proj}).
Since $G$ is connected,
$\overline{O}(\emptyset)={\cal B}$ is irreducible.
Since the reduced fibers of $\pi_1$ are $\cong \ppp^1$
(see Proposition \ref{lemma:fiber}),
it follows that $\overline{O}(s_1)$ is also irreducible,
by Lemma 5.3 in Debarre \cite{De}.
The claim follows by induction.
\qed

\begin{remark}
By induction, it follows that 
${\rm dim}(\overline{O}(s_1,\dots ,s_n))=n+{\rm dim}({\cal B})$. 
\label{remark:dim}
\end{remark} 

Next, we prove Lemma \ref{lemma:action}.

\begin{lemma}
\label{lemma:demazure}
Let $H_n=G\times P_1\times P_2\times \dots \times P_n$,
where $B$ is our fixed Borel, and
$P_i=\langle B,\, ^{s_i}B\rangle$, $1\leq i\leq n$.
Let $H_n$ operate\footnote{This is not an action.} 
on ${\cal B}^{n+1}$ by
\[
(g,p_1,...,p_n):(B_0,...,B_n)\longmapsto
(\, ^{g}B_0,\, ^{gp_1}B_1,\dots,\, ^{gp_1...p_n}B_n),
\]
for $p_i\in P_i$, $1\leq i\leq n$, and $g\in G$.
\begin{itemize}
\item[(a)]
 $\overline{O}(s_1,\dots ,s_{n})$ is the orbit of
$H_n$ operating on $x_{n+1}=(B,...,B)\in {\cal B}^{n+1}$:
$\overline{O}(s_1,\dots ,s_{n})=H_nx_{n+1}$.

\item[(b)]
There is a projection 
$\rho_n:H_n\rightarrow H_{n-1}$
such that, for all $h\in H_n$,
$\pi_n(h\cdot x_{n+1})=\rho_n(h)\cdot \pi_n(x_{n+1})=
\rho_n(h)\cdot x_{n}$.

\item[(c)]
Consider the projection $pr_1:\overline{O}(s_1,\dots ,s_{n})
\rightarrow {\cal B}$, sending
$(B_0,...,B_n)\longmapsto B_0$. Then
$pr_1^{-1}(B)=H'_nx_{n+1}$,
where $H'_n=P_1\times P_2\times \dots \times P_n$
and the operation of $H'_n$ is via the natural
embedding $H'_n\hookrightarrow H_n$.

\item[(d)]
Let $h,h'\in H'_n$ be given by 
$h=(p_1,...,p_n)$ and $h'=(p'_1,...,p'_n)$,
where $p_i,p'_i\in P_i$. We have
$h\cdot x_{n+1}=h'\cdot x_{n+1}$ if and only if
$p'_ip_i^{-1}\in \, ^{p_{i-1}...p_1}B$, for all
$1\leq i\leq n$ (when $i=1$, this should be interpreted as
$p'_1p_1^{-1}\in B$).

\item[(e)] 
Define a relation $\sim$ on $H'_n$ as follows: if
$h,h'\in H'_n$ are given by 
$h=(p_1,...,p_n)$ and $h'=(p'_1,...,p'_n)$,
where $p_i,p'_i\in P_i$, then define
$h\sim h'$ if and only if
$p'_ip_i^{-1}\in \, ^{p_{i-1}...p_1}B$, for all
$1\leq i\leq n$ (when $i=1$, this should be interpreted as
$p'_1p_1^{-1}\in B$).
This is an equivalence relation. 
\[
pr_1^{-1}(B)\cong H'_n/\sim,
\]
as sets.
\end{itemize}
\end{lemma}

\begin{remark}
(1) The proof below shows how to realize 
the Demazure varieties denoted $Z_K$ in \S 3 of \cite{D} 
as fibers in $\overline{O}(s_1,\dots ,s_n)$.

(2) If, in addition,
$s_1,s_1s_2,\dots ,s_1\dots s_n$ is a Hamiltonian path in the
Cayley graph of the Weyl group $W$ generated by the 
simple reflections of $B$ then this fact also 
explains how to embed Langlands' variety of stars
$S$ \cite{L}
into $\overline{O}(s_1,\dots ,s_n)$.

(3) Note that $Z_n$ (and, by Corollary \ref{cor:irred},
$\overline{O}(s_1,\dots ,s_n)$)
is irreducible, since $G$ is connected.
\end{remark}

\pf
(a):
Consider the $H_n$-orbit of $x_{n+1}$: $Z_n=H_nx_{n+1}$.
We claim
$Z_n=\overline{O}(s_1,\dots ,s_n)$.

By construction, $Z_n\subset \overline{O}(s_1,\dots ,s_n)$.
The projection 
$\pi_n':Z_n=H_nx_{n+1}\rightarrow Z_{n-1}=H_{n-1}x_{n}$
sends the point 
$(\, ^{g}B,\, ^{gp_1}B,\dots,\, ^{gp_1...p_n}B)\in Z_n$
to 
\newline
$(\, ^{g}B,\, ^{gp_1}B,\dots,\, ^{gp_1...p_{n-1}}B)\in Z_{n-1}$.

We identify the fiber of $\pi_n'$ 
over $x_n\in Z_{n-1}$, $(\pi_n')^{-1}(x_n)$.
By definition, 
\[
\begin{array}{c}
(\pi_n')^{-1}(x_n)=\\
\{(\, ^{g}B,\, ^{gp_1}B,\dots,\, ^{gp_1...p_n}B)\in Z_n
\ |\ ^{g}B=B,\, ^{gp_1}B=B,\dots,\, ^{gp_1...p_{n-1}}B=B\}.
\end{array}
\]
Since $B$ is it's own normalizer, 
these conditions become $g\in B$,
$gp_1\in B$ hence $p_1\in B$, ..., 
$gp_1...p_{n-1}\in B$ hence $p_{n-1}\in B$.
The only condition left unaffected then is that on $p_n$,
so $p_n\in P_n/B$.
We therefore see that this reduced 
fiber is $\cong P_n/B$, which is $\cong \ppp^1$
as in the proof of Proposition \ref{lemma:fiber}.
Also, note that the projection $\pi_n'$
is compatible with the inclusion
$Z_n\subset \overline{O}(s_1,\dots ,s_n)$.
Thus the fiber of $\pi_n'$ over $x_n$ is contained in the
corresponding fiber of $\pi_n$. 
As both of these fibers are
isomorphic to $\ppp^1$, they must be equal.
By the inclusion, each fiber of $\pi_n'$ is contained in the
corresponding fiber of $\pi_n$. We claim that
each fiber of $p'_n$ is of the form
$h\cdot (\pi'_n)^{-1}(x_n)$, for some $h\in H_n$.
Thus all fibers of $p'_n$ are isomorphic to $\ppp^1$.
Thus the fiber of $\pi_n'$ over any point is contained in the
corresponding fiber of $\pi_n$. 
As both of these fibers are
isomorphic to $\ppp^1$, they must be equal.

On the other hand, $Z_0\cong G/B\cong {\cal B}$.
Since the fibers of $\pi_1$ and $\pi_1'$ are equal, 
we have $Z_1=\overline{O}(s_1)$. By induction, we have
$Z_n=\overline{O}(s_1,\dots ,s_n)$. This proves (a).

(b), (c): An easy consequence of (a).

(d): Suppose 
\[
(\, B,\, ^{p_1}B,\dots,\, ^{p_1...p_n}B)
=
(\, B,\, ^{p'_1}B,\dots,\, ^{p'_1...p'_n}B).
\]
Since $B$ is its own normalizer in $G$,
$\, ^{p_1}B=\, ^{p'_1}B$ implies $p'_1p_1^{-1}\in B$.
Similarly, 
$\, ^{p_1p_2}B=\, ^{p'_1p'_2}B$ implies $p'_2p_2^{-1}\in \, ^{p_1}B$.
The claim in the lemma follows by an induction argument.

(e): To check that $\sim$ is an equivalence relation, note that
the condition 
\[
p_i(p'_i)^{-1}\in \, ^{p_{i-1}...p_1}B,\ \ \ \ \forall 1\leq i\leq n,
\]
is equivalent to 
\[
(p_1...p_i)(p'_1...p'_i)^{-1}\in B,\ \ \ \ \forall 1\leq i\leq n.
\]
Each of this last set of conditions is easily seen to be
both symmetric and transitive.

The last claim of (e) follows from (d) and the fact that $\sim$ is an 
equivalence relation.
\qed

\begin{proposition}
\label{prop:non-sing}
Assume that $k$ has characteristic $0$.
Each schematic fiber of $\pi_n$ is non-singular.
\end{proposition}

We conjecture this holds in characteristic $p>0$.

\pf
First, we prove the above proposition in the case $n=1$.
The property of being ``geometrically reduced''
holds on an open set (Theorem 12.2.4 in ch IV of \cite{EGA}).
We will show that, as a consequence of more general results
proven in section \S \ref{sec:appendix},
this open set contains the generic point, hence is dense
(and thus is non-empty).
As there is a transitive group action of $G$ on
the set of these fibers (Lemma \ref{lemma:demazure} (b)), 
all the fibers are isomorphic (as schemes), hence 
they must all be geometrically reduced. 
We conclude from
Proposition \ref{lemma:fiber} that these fibers are all 
(reduced schemes) isomorphic to $\ppp^1$.

Now the proposition follows from the schematic analog of
(\ref{eqn:fibers}), 
which is a special case of
(\ref{eqn:fibers2}).
\qed

\section{The smoothness result}
\label{sec:3}

We shall use the following result.

\begin{lemma}
\label{lemma:proj}
For each $(w_1,...,w_n)\in W^n$, we have
\[
\overline{O}(w_1,\dots ,w_n)=\overline{O}(w_1,\dots ,w_{n-1})
\times_{\cal B} \overline{O}(w_n),
\]
where the morphism
$\overline{O}(w_1,\dots ,w_{n-1})\rightarrow {\cal B}$ is given by
$(B_0,\dots ,B_{n-1})\longmapsto B_{n-1}$ and
the morphism
$\overline{O}(w_n)\rightarrow {\cal B}$ is given by
$\pi_1:(B'_0,B'_1)\longmapsto B'_0$.
The projection morphism $\pi_n$ is flat.
\end{lemma}

\pf
Let $X'=\overline{O}(w_1,\dots ,w_n)$,
$X=\overline{O}(w_1,\dots ,w_{n-1})$,
$Y'=\overline{O}(w_n)$,
and let $Y={\cal B}$. 
The fact that $X'=X\times_Y Y'$ follows from the definition 
of $X'=\overline{O}(w_1,\dots ,w_n)$
and $\times_Y$.
We must show that $\pi_1:Y'\rightarrow Y$ is flat.

A well-known theorem implies 
that $\pi_1$ is generically flat (see Theorem 6.9.1 in ch IV of \cite{EGA} or
for example, Corollary III.10.7 in \cite{H}). Since $G$ acts transitively on 
$Y$ and $\pi_1$ is $G$-equivariant
by Lemma \ref{lemma:demazure} (b), $\pi_1$ must be flat
over $O(w_n)\subset Y'$.

Note that $G$ acts transtively on the fibers of $\pi_1$.
In particular, given any two $x,x'\in \overline{O}(w_n)$
with, say $x\in \pi_1^{-1}(z)$ and $x'\in \pi_1^{-1}(z')$
for $z,z'\in {\cal B}$, there is a $g\in G$ such that
$gz=z'$. In other words, we can ``translate'' any fiber
of any point in $\overline{O}(w_n)$
to the fiber of a point in $O(w_n)$. 
Since $\pi_1$ is flat over $O(w_n)$, 
it follows $\pi_1$ is flat everywhere on $\overline{O}(w_n)$.

The conclusion that $\pi_n$ is flat follows
from ``base change'' (Proposition III.9.2(b) in Hartshorne
\cite{H}).
\qed

\begin{remark}
To help visalize the geometry,
we construct an open $U\subset Y'$ containing $Y'-O(w_n)$.
Let $B'$ denote a Borel opposite to $B$, let $N'$ denote the
unipotent radical of $B'$, and let 
$P_{w}=\langle B,\, ^{w}B\rangle$ ($w\in W$). Define
\begin{equation}
\label{eqn:U}
U(w)=\{(^nB,\, ^{np}B)\ |\ n\in N',\ p\in P_{w}\}.
\end{equation}
The map $N'\times P_{w_n}/B\rightarrow U(w_n)$,
defined by $(n,p)\longmapsto (^nB,\, ^{np}B)$,
is an isomorphism. 
\end{remark}

Next we prove Proposition \ref{prop:smooth}:
for ${\rm char}(k)=0$,
the closed variety of Borels of relative position
$\vec{s}$ is non-singular and $\pi_n$ is smooth.

\pf We prove this by induction on $n$,
using Lemma \ref{lemma:1.1}.
 
Let $n=1$. The reduced fibers are all $\ppp ^1$ by 
Proposition \ref{lemma:fiber}. 
Since $G$ acts transitively on ${O(s_1)}$,
it must be non-singular. Each point in $U(s_1)$ is
smooth, where $U$ is as in (\ref{eqn:U}),
so $\overline{O}(s_1)$ is non-singular.
Each schematic fiber of $\pi_1$ is non-singular, by Proposition
\ref{prop:non-sing}.
Since each schematic fiber is non-singular, and its
reduced scheme is $\ppp^1$, they must all be (isomorphic to) $\ppp^1$. 

Since the schematic fibers are non-singular and
$\pi_1$ is flat, by Theorem 3' in \S III.10 of Mumford \cite{Mum},
$\pi_1$ is smooth.
 
Now assume $n>1$. By the induction hypothesis,  
$\overline{O}(s_1,\dots ,s_{n-1})$ 
is non-singular. As above, it suffices to check that 
$\pi_n$ is flat (by Mumford's result and
Proposition \ref{prop:non-sing}). So by Lemma \ref{lemma:proj}
again,
$\pi_n$ is smooth. By the induction hypothesis and Lemma
\ref{lemma:1.1}, $\overline{O}(s_1,\dots ,s_{n})$ 
is non-singular.

The schematic fibers are all $\ppp ^1$ by 
the same argument as in the case $n=1$.

\qed 
 
\section{Remarks on $\ppp^1$-bundles}
\label{sec:4}

In this section, we discuss Question \ref{thrm:main}.

First, we recall the definition of a $\ppp^1$-bundle.
A {\bf $\ppp^1$-bundle} $E$ over $X$ is 
(1) a morphism $\varphi :E\rightarrow X$,
(2) an open covering $\{U_\alpha\}$ of $X$ 
and a collection of isomorphisms
$\varphi_\alpha:\varphi^{-1}(U_\alpha)\rightarrow 
U_\alpha\times\ppp^1$ such that
(a) the following diagram commutes
\[
\begin{CD}
\varphi^{-1}(U_\alpha) @>\varphi_\alpha >>U_\alpha\times\ppp^1 \\
\varphi @VVV @VV pr_1 V\\
U_\alpha @= U_\alpha,
\end{CD}
\]
and (b) for all $\alpha, \beta$,
\[
\varphi_\beta\varphi_\alpha^{-1}:
U_\alpha\cap U_\beta\times\ppp^1\rightarrow
U_\alpha\cap U_\beta\times\ppp^1
\]
may be regarded as an element of 
$\Gamma(U_\alpha\cap U_\beta,PGL(2))$.
This definition is essentially a special case
of that in Hartshorne \cite{H}, page 170.

Regarding the question of whether or not
\[
\overline{O}(s_1,\dots ,s_n)\stackrel{\pi_n}{\rightarrow}
\overline{O}(s_1,\dots ,s_{n-1})
\]
is a $\ppp^1$-bundle (so $\varphi=\pi_n$ in part (1) of the definition above),
let us consider the case where $n=1$:
$\overline{O}(s)
\stackrel{\pi_1}{\rightarrow}
\overline{O}(\emptyset)={\cal B}$.
Write $(B_0,B_1)\in \overline{O}(s)$ as
$B_0=gB$, $B_1=gsB$, for some simple reflection $s$
(with respect to the Borel $B$) and some $g\in G$.
Under $\pi_1$, the pair $(B_0,B_1)\in \overline{O}(s)$ is sent to
$B_0\in \overline{O}(\emptyset)$. 

The map $N'\times P_{s}/B\rightarrow U(s)$,
defined by $(n,p)\longmapsto (^nB,\, ^{np}B)$,
is an isomorphism, where $U(s)$ 
is as in (\ref{eqn:U}) and $P_{s}=\langle B,\, ^{s}B\rangle$.
Let $\psi_n:\, ^nP_{s}/\, ^nB\rightarrow \ppp^1$ denote an isomorphism.
Consider the map $\phi:N'\times \ppp^1\rightarrow \overline{O}(s)$
defined by
\[
\begin{array}{ccc}
N'\times \ppp^1\rightarrow & N'\times\, 
^nP_s/\, ^nB
\rightarrow & \overline{O}(s)\\
(n,x)\longmapsto & (n, \psi_n(x))\longmapsto &(^{n}B,\, ^{n\psi(x)}B)
\end{array}
\]
In this case, the next step to addressing the question 
is to verify that the following diagram commutes:
\[
\begin{CD}
N'\times \ppp^1 @>\phi >> \pi_1^{-1}(N')\\
pr_1@VVV @VV\pi_1 V\\
N' @= N'.
\end{CD}
\]
This follows from the definition.
Therefore, we have some evidence\footnote{
We have not even shown, among other things, that
the above map $\phi$ represents a morphism.} 
that condition (2a) of the definition above holds,
in the case $n=1$.

Although our discussion above does not settle the question,
from what we have shown it follows that 
$\overline{O}(s)$ is a uniruled variety, in the sense of
Debarre \cite{De}. 

\section{Deligne-Lusztig varieties}
\label{sec:5}

We prove Corollary \ref{cor:DL},
under the assumption that Proposition \ref{prop:smooth} 
holds for $k=\overline{\fff_q}$.

Let $\vec{s}=(s_1,\dots ,s_n)\in S^n$ be simple reflections 
with respect to
a fixed Borel $B$ in a connected semisimple group $G$ defined over
a finite field $\fff_q$. Let 
$\overline{X}(\vec{s})$ be the Deligne-Lusztig variety,
as in \S \ref{sec:1}.
From \cite{DL}, page 151, we know that this variety is
given by the fiber product
\[
\overline{X}(\vec{s})=
\overline{O}(\vec{s})\times_{{\cal B}\times {\cal B}}\Gamma,
\]
where $\Gamma$ is the graph of the Frobenius morphism $F$. 
By Lemma \ref{lemma:1.1}, Proposition \ref{prop:smooth}, and
the fact that the intersection is transversal
(see \cite{DL}, Lemma 9.11), 
$\overline{X}(\vec{s})$ is non-singular.
\qed

\section{Appendix: Technical background}
\label{sec:appendix}

This section is devoted to finishing the proof of Proposition
\ref{prop:non-sing}.

First, we record a simple observation used repeatedly in the
previous sections.

\begin{lemma} Let $f:X\rightarrow Y$ 
be a smooth morphism of 
varieties defined over $k$. If $Y$ is non-singular 
then $X$ is also non-singular. 
\label{lemma:1.1} 
\end{lemma} 

This is well-known. 

\subsection{Standard properties of fibers}
\label{subsec:fibers}

We recall some standard facts due to Grothendieck \cite{EGA}.

From EGA I, 3.6.3:
Assume that $X$ and $Y$ are affine schemes, that $F:X\rightarrow Y$ is a morphism,
and that 
$U\subset X$ is open. The following diagram is commutative.
\[
\begin{CD}
X @<<< U\\
f@VVV @VV{f|_U}V\\
Y @= Y,
\end{CD}
\]
where the top arrow is inclusion. For $y\in f(U)$, the morphism  
\[
f^{-1}(U)\cap U \rightarrow (f|U)^{-1}(y)
\]
induced from $f$, is an isomorphism.

From EGA I, 3.6.4:
Let $Y'$ be an irreducible component of $Y$
and let $f:X\rightarrow Y$. Let
\[
X_{(Y')}=X\times_Y Y'\stackrel{f_{(Y')}}{\rightarrow} Y'.
\]
We have $f^{-1}(Y')=X'=X_{(Y')} $. Moreover, 
\[
\begin{CD}
X_{(Y')} @>>> X\\
@Vf_{(Y')}VV @VVfV \\
Y' @>j>> Y,
\end{CD}
\]
where $j$ is inclusion. Then
\begin{equation}
f^{-1}(y)\leftarrow f_{(Y')}^{-1}(y'),
\label{eqn:fibers2}
\end{equation}
is an isomorphism as schemes
(``transitivity of fibers'').

\subsection{Main lemma}

Here is our main lemma.

\begin{lemma}
\label{lemma:main}
Let $D$ be an integral domain with field of fractions,
$K=Frac(D)$. Let $A$ be a $D$-algebra, flat over $D$
and let $A_K=A\otimes_D K$. Then
\begin{itemize}
\item[(i)]
if $A$ is an integral domain then so is $A_K$,
\item[(ii)]
if $A$ is reduced then so is $A_K$.
\end{itemize}

\end{lemma}

\begin{remark}
Here is an example to show that the extension of
(b) to ``if $A$ is reduced then $A_K$ is geometrically 
reduced'' is false.
(This counterexample suggests that the ``${\rm char}(k)=0$'' hypothesis 
in Proposition \ref{prop:perfect} below cannot be removed,
at least not without adding some other condition.)
Let $k=\overline{\fff_2}$, $D=k[x]$,
$K=k(x)$, so $A_K=k(x)[y]/(y^2-x)$.
Consider the purely inseparable extension $K=k(\sqrt{x})$
of $K=k(x)$. Then $A_{K'}=k(\sqrt{x})[y]/(y^2-x)$.Consider the 
non-zero element $a=y+\sqrt{x}\in A_{K'}$. 
This element satisfies $a^2=0$, so $A_K$ is not 
geometrically reduced.
\end{remark}

\pf
By flatness, the map
\[
A=A\otimes_D D\rightarrow A_K
\]
is injective, so we may view $A$ as a subset of $A_K$.

For (i): Assume $A_K$ is not a domain and assume,
that $(\sum x\otimes y)\cdot (\sum x'\otimes y')=0$, for some
non-zero $\sum x\otimes y, \sum x'\otimes y'\in A_K$.
Eliminating denominators, we may assume $y,y'\in D$.
This gives divisors of zero in $A\otimes_D D=A\subset A_K$.
This proves (i).

The proof of (ii) is similar, hence omitted.

\qed

\begin{proposition}
\label{prop:perfect}
Assume that $k$ has characteristic $0$.
Let $X$ and $Y$ be reduced schemes defined $k$ and let
$f:X\rightarrow Y$ be a proper, surjective, flat morphism of finite
presentation.
Then
\begin{itemize}
\item[(i)]
The generic fiber of $f$ is geometrically reduced.

\item[(ii)]
The general fiber of $f$ is geometrically reduced (i.e., all fibers over a 
dense open set in $Y$ are reduced).

\end{itemize}
\end{proposition}

\begin{remark}
 We only need to assume $f$ is flat for part (i).
\end{remark}

\pf 
We only prove the generic fiber of $f$ is reduced here.
The next subsection extends the proof to the
geometrically reduced situation.

Part (ii) follows from (i) and \cite{EGA} IV, \S 12.2.4.

It remains to prove (i). Using the standard properties of fibers 
(in \S \ref{subsec:fibers})
and the fact that $Y$ is Noetherian, we reduce to the case
where $X$ and $Y$ are affine and $Y$ is integral. 
Appealing to Lemma \ref{lemma:main} proves (i)
with ``geometrically reduced''
replaced by ``reduced''.

It remains to show that the generic fiber of $f$ is
geometrically reduced. By another result of Grothendieck,
it suffices to show the fiber remains reduced after passing to a
finite purely inseparable extension. 
But if $k$ is perfect then every finite purely inseparable extension
is trivial.
\qed

\end{document}